\pdfoutput=1
\RequirePackage{ifpdf}
\ifpdf % We are running pdfTeX in pdf mode
\documentclass[pdftex]{sigma}
\else
\documentclass{sigma}
\fi
\usepackage[all]{xy}

\newcommand{\C}{\mathbb{C}}
\newcommand{\R}{\mathbb{R}}

\newcommand{\Z}{\mathbb{Z}}

\newcommand{\DD}{\mathbb{D}}

\newcommand{\PP}{\mathbb{P}}
\newcommand{\U}{\mathrm{U}}

\numberwithin{equation}{section}

\newtheorem{Theorem}{Theorem}[section]
\newtheorem{Proposition}[Theorem]{Proposition}
\newtheorem{Conjecture}[Theorem]{Conjecture}
 { \theoremstyle{definition}
\newtheorem{Definition}[Theorem]{Definition}
\newtheorem{Example}[Theorem]{Example}
\newtheorem{Remark}[Theorem]{Remark} }

\begin{document}

%\allowdisplaybreaks

\newcommand{\arXivNumber}{1612.04623}

\renewcommand{\thefootnote}{}

\renewcommand{\PaperNumber}{024}

\FirstPageHeading

\ShortArticleName{Doran--Harder--Thompson Conjecture via SYZ Mirror Symmetry: Elliptic Curves}

\ArticleName{Doran--Harder--Thompson Conjecture\\ via SYZ Mirror Symmetry: Elliptic Curves\footnote{This paper is a~contribution to the Special Issue on Modular Forms and String Theory in honor of Noriko Yui. The full collection is available at \href{http://www.emis.de/journals/SIGMA/modular-forms.html}{http://www.emis.de/journals/SIGMA/modular-forms.html}}}

\Author{Atsushi KANAZAWA}

\AuthorNameForHeading{A.~Kanazawa}

\Address{Department of Mathematics, Kyoto University,\\ Kitashirakawa-Oiwake, Sakyo, Kyoto, 606-8502, Japan}
\Email{\href{akanazawa@math.kyoto-u.ac.jp}{akanazawa@math.kyoto-u.ac.jp}}

\ArticleDates{Received December 20, 2016, in f\/inal form April 06, 2017; Published online April 11, 2017}

\Abstract{We prove the Doran--Harder--Thompson conjecture in the case of elliptic curves by using ideas from SYZ mirror symmetry. The conjecture claims that when a Calabi--Yau manifold $X$ degenerates to a union of two quasi-Fano manifolds (Tyurin degeneration), a~mirror Calabi--Yau manifold of $X$ can be constructed by gluing the two mirror Landau--Ginzburg models of the quasi-Fano manifolds. The two crucial ideas in our proof are to obtain a complex structure by gluing the underlying af\/f\/ine manifolds and to construct the theta functions from the Landau--Ginzburg superpotentials. }


\Classification{53D37; 14J33; 14J32; 14J45; 14D06}

\renewcommand{\thefootnote}{\arabic{footnote}}
\setcounter{footnote}{0}

\section{Introduction}
The aim of this short article is to prove the Doran--Harder--Thompson conjecture in the case of elliptic curves by using ideas from SYZ mirror symmetry.

Given a Tyurin degeneration of a Calabi--Yau manifold $X$ to a union of two quasi-Fano mani\-folds $X_1 \cup_Z X_2$ intersecting along a common smooth anti-canonical divisor $Z \in |{-}K_{X_i}|$ for $i=1,2$, it is natural to investigate a potential relationship between geometry of the Calabi--Yau manifold $X$ and that of the quasi-Fano manifolds $X_1$ and $X_2$. Motivated by the works of Dolgachev~\cite{Dol}, Tyurin~\cite{Tyu} and Auroux~\cite{Aur2}, recently Doran--Harder--Thompson proposed a~remarkable conjecture (Conjecture~\ref{DHT}), which builds a bridge between mirror symmetry for the Calabi--Yau manifold $X$ and that for the quasi-Fano manifolds~$X_1$ and~$X_2$~\cite{DHT}. It claims that we should be able to glue together the mirror Landau--Ginzburg models $W_i\colon Y_i \rightarrow \C$ of the pair $(X_i,Z)$ for $i=1,2$ to construct a mirror Calabi--Yau manifold $Y$ of $X$ equipped with a~f\/ibration $W\colon Y\rightarrow \PP^1$. They provided supporting evidence for the conjecture in various dif\/ferent settings. For instance it was shown that under suitable assumptions we can glue together the Landau--Ginzburg models $W_i\colon Y_i \rightarrow \C$ for $i=1,2$ to obtain a $C^\infty$-manifold $Y$ with the expected topological Euler number $\chi(Y)=(-1)^{\dim X}\chi(X)$. Thus the topological version of the conjecture is essentially proven. However, the real dif\/f\/iculty of this conjecture lies in constructing $Y$ as a~complex manifold, which should be mirror to the symplectic manifold $X$ (or vice versa).

In this article, we will prove the conjecture in the case of elliptic curves, beginning with a~symplectic elliptic curve $X$ and constructing the mirror complex elliptic curve $Y$. In order to obtain the correct complex manifold $Y$ by gluing the mirror Landau--Ginzburg models \smash{$W_i\colon Y_i\rightarrow \C$} for $i=1,2$, we need to keep track of the subtle complex structures of the (not necessarily algebraic) K\"ahler manifold $Y_i$ for $i=1,2$. To this end, we f\/ind ideas of SYZ mirror symmetry very useful.

The key idea in our proof is twofold. The f\/irst is to obtain the correct complex structure by gluing the underlying af\/f\/ine base manifolds of $X_1$ and $X_2$ in SYZ mirror symmetry. This is based on the philosophy that a Tyurin degeneration of a Calabi--Yau manifold~$X$ can be thought to be f\/ibred over a Heegaard splitting of the base $B$ of a special Lagrangian torus f\/ibration $\phi\colon X\rightarrow B$. The second is to construct theta functions out of the Landau--Ginzburg superpotentials. As a~corollary, we observe that the product formulae of the theta functions are the manifestation of quantum corrections appearing in SYZ mirror symmetry.

It is worth mentioning that a variant of the conjecture was discussed in the work of Auroux~\cite{Aur2}. He studied a $\Z/2\Z$-invariant version of mirror symmetry (Calabi--Yau double covers) and investigated elliptic curves from a dif\/ferent perspective from ours \cite[Example~3.2]{Aur2}. Advantages of the present work are f\/irstly to work with not necessarily identical af\/f\/ine manifolds (doubling) and secondly to construct theta functions out of the Landau--Ginzburg superpotentials in an interesting way based on the geometry of the conjecture.

\subsection*{Structure of article}
We will provide a self-contained description of the subjects for completeness. Section \ref{Tyurin} builds basic setup and formulates the main conjecture. Section~\ref{SYZ} reviews basics of SYZ mirror symmetry both in the Calabi--Yau and quasi-Fano settings. Section \ref{Gluing} is the main part of this article and proves the conjecture in the case of elliptic curves. Section~\ref{Direction} comments on further research directions.

\section{Doran--Harder--Thompson conjecture} \label{Tyurin}
In this section, we will provide background materials and review the Doran--Harder--Thompson conjecture, following the original article~\cite{DHT}. In general mirror symmetry is a conjecture about a Calabi--Yau manifold near a large complex structure limit, which is thought to be a maximal degeneration, in the complex moduli space. However in this article we will be interested in another class of loci in the complex moduli space, where a Calabi--Yau manifold degenerates to a~union of two quasi-Fano manifolds.

\subsection{Tyurin degeneration}
A Calabi--Yau manifold $X$ is a compact K\"ahler manifold such that the canonical bundle is trivial $K_X =0$ and $H^i(X,\mathcal{O}_X)=0$ for $0 < i <\dim X$. A quasi-Fano manifold $X$ is a smooth variety $X$ such that $|{-}K_X|$ contains a smooth Calabi--Yau member and $H^i(X,\mathcal{O}_X)=0$ for $0<i$. A~Tyurin degeneration is a degeneration $\mathcal{X}\rightarrow \Delta$ of Calabi--Yau manifolds over the unit disc $ \Delta=\{|z|<1\} \subset \C$, such that the total space $\mathcal{X}$ is smooth and the central f\/ibre $\mathcal{X}_0=X_1 \cup_Z X_2$ is a~union of two quasi-Fano manifolds $X_1$ and~$X_2$ intersecting normally along a~common anti-canonical divisor $Z \in |{-}K_{X_i}|$ for $i=1,2$. Conversely, we have the following result of Kawamata--Namikawa \cite{KN}, which is slightly modif\/ied for our setting.

\begin{Theorem}[{\cite[Theorem~4.2]{KN}}] \label{KN}
Let $X_1$ and $X_2$ be quasi-Fano manifolds and $Z \in |{-}K_{X_i}|$ a~common smooth anti-canonical divisor for $i=1,2$. Assume that there exist ample divisors $D_i \in \operatorname{Pic}(X_i)$ for $i=1,2$ which restrict to an ample divisor $D_1|_Z=D_2|_Z$ on~$Z$. Then the union $X_1 \cup_Z X_2$ of $X_1$ and~$X_2$ intersecting normally along~$Z$ is smoothable to a~Calabi--Yau manifold~$X$ if and only if $N_{Z/X_1}\cong N_{Z/X_2}^{-1}$ $(d$-semistability$)$. Moreover the resulting Calabi--Yau manifold~$X$ is unique up to deformation.
\end{Theorem}

A Tyurin degeneration is thought to be a complex analogue of a Heegaard splitting of a~compact oriented real 3-fold without boundary. Based on this analogy, in his posthumous article \cite{Tyu} Tyurin proposed to study geometry of a Calabi--Yau 3-fold by using that of quasi-Fano 3-folds when they are related by a Tyurin degeneration.

\subsection{Mirror symmetry for quasi-Fano manifolds}
We consider a pair $(X,Z)$ consisting of a quasi-Fano manifold $X$ and an anti-canonical divisor $Z\in |{-}K_X|$. The complement $X\setminus Z$ can be thought of as a log Calabi--Yau manifold as there exists a nowhere vanishing volume form $\Omega$ on $X\setminus Z$ with poles along~$Z$.

\begin{Example} \label{toric}
Let $X$ be a toric Fano $n$-fold and $Z$ the toric boundary, which is the complement of the dense torus $(\mathbb{C}^\times)^n \subset X$. Then $X\setminus Z=(\mathbb{C}^\times)^n$ carries a standard holomorphic volume form $\Omega=\wedge_{i=1}^{n}\sqrt{-1}d\log z_i$, where $(z_i)$ are the coordinates of $(\mathbb{C}^\times)^n$.
\end{Example}

\begin{Definition} A Landau--Ginzburg model is a pair $(Y,W)$ of a K\"ahler manifold $Y$ and a~holomorphic function $W\colon Y\rightarrow \mathbb{C}$, which is called a superpotential.
\end{Definition}

It is classically known that there is a version of mirror symmetry for Fano manifolds together with an anti-canonical divisor. We expect that such mirror symmetry should hold also for quasi-Fano manifolds (or even for varieties with ef\/fective anti-canonical divisors~\cite{Aur1}). Here we formulate a mirror conjecture for quasi-Fano manifolds (see for example Katzarkov--Kontsevich--Pantev~\cite{KKP}, Harder~\cite{Har}).

\begin{Conjecture} \label{LG MS}For a pair $(X,Z)$ of a quasi-Fano $n$-fold $X$ and a smooth anti-canonical divisor $Z \in |{-}K_X|$, there exists a Landau--Ginzburg model $(Y,W)$ such that
\begin{enumerate}\itemsep=0pt
\item[$1)$] $\sum_{j}h^{n-i+j,j}(X)=h^i(Y,W^{-1}(s))$ for a regular value $s\in \C$ of $W$,
\item[$2)$] the generic fibres of $W$ and the generic anti-canonical hypersurfaces in $X$ are mirror families of compact Calabi--Yau $(n-1)$-folds,
\end{enumerate}
where $h^i(Y,W^{-1}(s))$ is the rank of the relative cohomology group $H^i(Y,W^{-1}(s))$. The pair $(Y,W)$ is called a mirror Landau--Ginzburg model of $(X,Z)$.
\end{Conjecture}

The anti-canonical divisor $Z$ can be thought of as an obstruction for the quasi-Fano mani\-fold~$X$ to be a Calabi--Yau manifold and~$W$ is an obstruction (or potential function) for the Floer homology of a Lagrangian torus in~$X$ to be def\/ined in the sense of Fukaya--Oh--Ohta--Ono~\cite{CO, FOOO} as we will see in the next section.

\begin{Example} \label{HVMirror}
Let $X=\mathbb{P}^1$ and $Z=\{0,\infty\}$ equipped with a toric K\"ahler form $\omega$. Then the mirror Landau--Ginzburg model of $(X,Z)$ is given by $(\mathbb{C}^\times, W=z +\frac{q}{z})$, where $q=\exp\big({-}\int_{\mathbb{P}^1}\omega\big)$.
One justif\/ication of this mirror duality is given by the ring isomorphism
\begin{gather*}
\operatorname{QH}\big(\mathbb{P}^1\big)=\mathbb{C}[H]/\big(H^2-q\big) \cong \mathbb{C}\big[z^{\pm1}\big]/\big(z^2-q\big) =\operatorname{Jac}(W).
\end{gather*}
Here $\operatorname{QH}(\mathbb{P}^1)$ is the quantum cohomology ring of $\mathbb{P}^1$ and $\operatorname{Jac}(W)$ is the Jacobian ring of the superpotential~$W$.
\end{Example}

Conjecture \ref{LG MS} can be generalized to the case when $Z$ is mildly singular, but the superpotential~$W$ will no longer be proper. In fact, it is anticipated that $Z$ is smooth if and only if~$W$ is proper.

\begin{Example} \label{P^2}
For $X=\mathbb{P}^2$ with a toric K\"ahler form $\omega$, we def\/ine $Z_0$ to be the toric boundary, $Z_1$~the union of a~smooth conic and a line intersecting 2 points, and $Z_2$ a nodal cubic curve, and~$Z_3$ a~smooth cubic curve. The mirror Landau--Ginzburg model of $(X,Z_0)$ is given by
\begin{gather*}
\left(Y_0=(\mathbb{C}^\times)^2, W_0=x+y +\frac{q}{xy}\right),
\end{gather*}
where $q=\exp\big({-}\int_{H}\omega\big)$ for the line class $H$. A generic f\/iber of $W_0$ is an elliptic curve with 3~punctures. On the other hand, the mirror Landau--Ginzburg model $(Y_i,W_i)$ of the pair $(X,Z_i)$ is a f\/iberwise partial compactif\/ication of $W_0$ such that a generic f\/iber of $W_i$ is an elliptic curve with $3-i$ punctures for $1 \le i \le 3$.
\end{Example}

\subsection{Doran--Harder--Thompson conjecture} \label{DHT}
Mirror symmetry for Calabi--Yau manifolds and that for (quasi-)Fano manifolds have been \mbox{studied} for a long time, but somewhat independently.
A natural question to ask is, {\it how are mirror symmetry for these manifolds related to each other?} Motivated by works of Dolgachev~\cite{Dol}, Tyurin~\cite{Tyu}, and Auroux~\cite{Aur2}, Doran--Harder--Thompson proposed the following remarkable conjecture, which we call the DHT conjecture for short.

\begin{Conjecture}[Doran--Harder--Thompson \cite{DHT}]
Given a Tyurin degeneration of a Calabi--Yau manifold $X$ to the union $X_1 \cup_Z X_2$ of quasi-Fano manifolds intersecting along their common smooth anti-canonical divisor~$Z$, then the mirror Landau--Ginzburg models $W_i\colon Y_i\rightarrow \C$ of $(X_i,Z)$ for $i=1,2$ can be glued together to be a Calabi--Yau manifold $Y$ equipped with a Calabi--Yau fibration $W\colon Y\rightarrow \PP^1$. Moreover, $Y$ is mirror symmetric to $X$.
\end{Conjecture}

\begin{figure}[htbp]\centering
 \includegraphics[width=50mm]{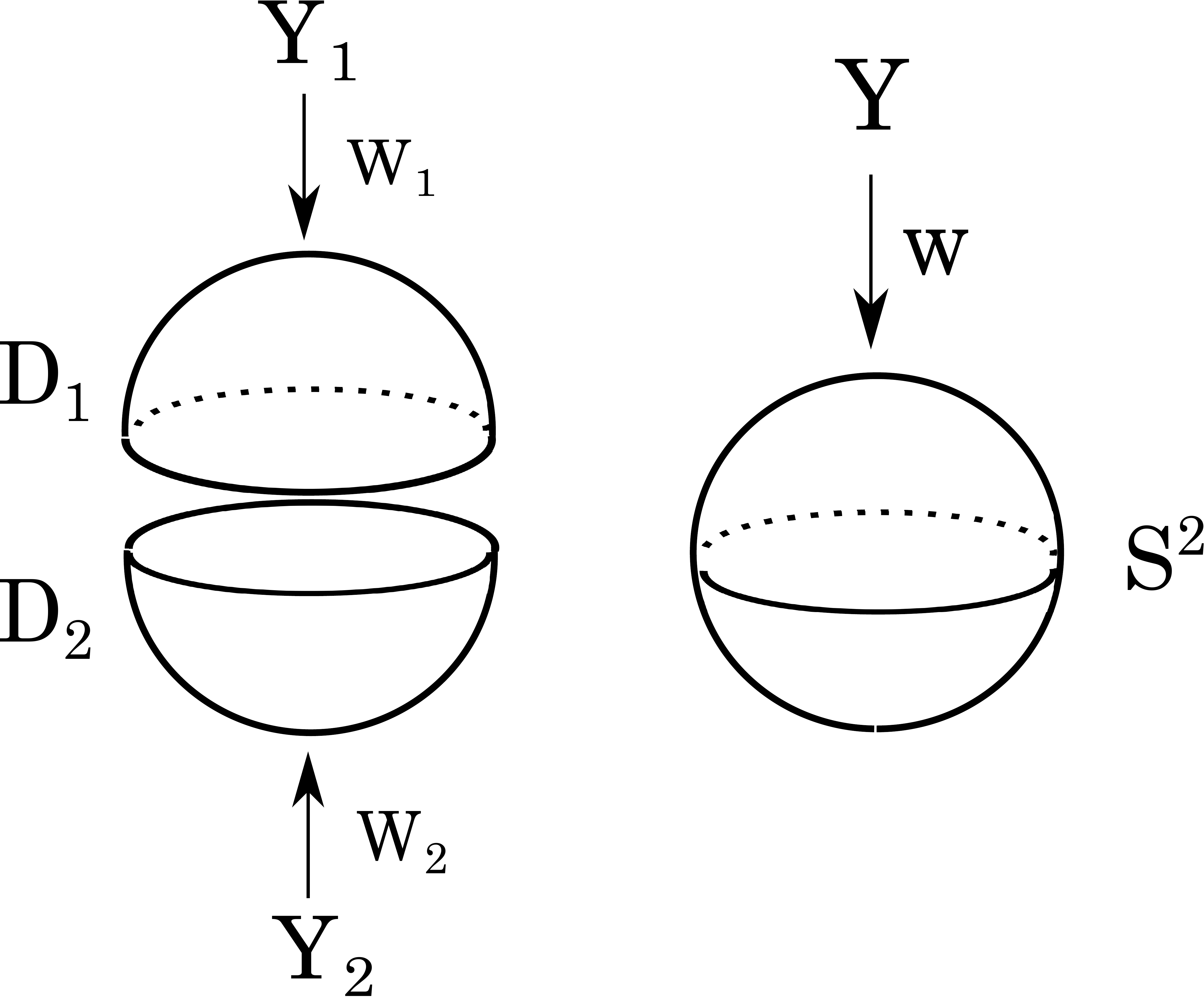}

 \caption{Gluing Landau--Ginzburg models.}\label{fig:gluing}
\end{figure}

The above gluing process can be understood as follows. We denote by $n$ the dimension of $X$ and by $Z_i^\vee$ a f\/iber of the superpotential $W_i$ mirror to a Calabi--Yau $(n-1)$-fold $Z$.
\begin{enumerate}\itemsep=0pt
\item Firstly, we assume that all the important information about the Landau--Ginzburg model $W_i\colon Y_i\rightarrow \C$ is contained in the critical locus of the superpotential $W_i$. Therefore, without much loss of information, we may replace it with a new Landau--Ginzburg model \smash{$W_i\colon Y_i\rightarrow \DD_i$} for a suf\/f\/iciently large disc $\DD_i$ which contains all the critical values.
\item Secondly, the Calabi--Yau manifolds $Z_1^\vee$ and $Z_2^\vee$ are both mirror symmetric to $Z$, and thus we expect that they are topologically identif\/ied\footnote{Two Calabi--Yau manifolds may be topologically dif\/ferent even if they share the same mirror manifold. There is no problem if $\dim X=1,2$ or $3$ for example.}.
\item Thirdly, Theorem \ref{KN} implies that we have $N_{Z/X_1}\cong N_{Z/X_2}^{-1}$ because $X_1 \cup_Z X_2$ is smoothable to a Calabi--Yau manifold $X$.
According to Kontsevich's homological mirror symmetry \cite{Kon}, we have an equivalence of triangulated categories
\begin{gather*}
\mathrm{D^bCoh}(Z)\cong \mathrm{D^bFuk}\big(Z_i^\vee\big).
\end{gather*}
Then the monodromy symplectomorphism on $Z_i^\vee$ associated to the anti-clockwise loop $\partial \DD_i$ can be identif\/ied
with the autoequivalence $(-)\otimes \omega_{X_i}[n]|_Z\cong (-)\otimes N_{Z/X_i}^{-1}[n]$ on $ \mathrm{D^bCoh}(Z)$ (see \cite{KKP, Sei} for details).
Therefore $N_{Z/X_1}\cong N_{Z/X_2}^{-1}$ implies, under mirror symmetry, that the monodromy action on $Z_1^\vee$
along the anti-clockwise loop $\partial \DD_1$ and that on $Z_2^\vee$ along the clockwise loop $-\partial \DD_2$ can be identif\/ied.
\end{enumerate}
Therefore, assuming various mirror symmetry statements, we are able to glue the f\/ibrations $W_i\colon Y_i\rightarrow \DD_i$ for $i=1,2$ along open neighborhoods of the boundaries $\partial \DD_1$ and $\partial \DD_2$ to construct a~$C^\infty$-manifold $Y$ equipped with a~f\/ibration $W\colon Y\rightarrow S^2$ (Fig.~\ref{fig:gluing}). Note that the smoothness of~$Z$ implies the compactness of $Y$, which follows from the properness of the superpotentials $W_i$ for $i=1,2$ (cf.\ Example~\ref{P^2}). The highly non-trivial part of the conjecture is that there exist a Calabi--Yau structure on~$Y$ and a complex structure on~$S^2$ in such a way that $W\colon Y\rightarrow \PP^1$ is holomorphic and $Y$ is mirror symmetric to the Calabi--Yau manifold~$X$. Moreover, a f\/iber of~$W$ is mirror to $Z$ by the above construction.

In \cite{DHT} the authors provide supporting evidence for this conjecture in various dif\/ferent settings, including Batyrev--Borisov mirror symmetry and Dolgachev--Nikulin mirror symmetry. For example, in the 3-fold case, the Hodge number mirror symmetry is shown to be equivalent to a version of Dolgachev--Nikulin mirror symmetry for K3 surfaces, provided that $Y$ admits a Calabi--Yau structure. Another important result is that, under reasonable assumptions, the resulting $C^\infty$-manifold $Y$ has the expected Euler number $\chi(Y)=(-1)^{\dim X}\chi(X)$. Thus the conjecture is essentially proven at the topological level. However, the real dif\/f\/iculty of the conjecture lies in constructing~$Y$ as a complex manifold, which should be mirror to the symplectic mani\-fold~$X$ (or vice versa). The goal of this article is to overcome this dif\/f\/iculty in the 1-dimensional case.

One warning is in order. The DHT conjecture is likely to be false unless we impose a~condition on how the Tyurin degeneration of a Calabi--Yau manifold $X$ occurs in the complex moduli space. For example, if the complex moduli space of $X$ is $1$-dimensional, the K\"ahler moduli space of a mirror Calabi--Yau manifold $Y$ is also $1$-dimensional and thus $Y$ cannot have a f\/ibration structure (unless it is 1-dimensional). The conjecture should be modif\/ied so that the Tyurin degeneration occurs in a locus which contains a large complex structure limit. In such a~case, the Calabi--Yau manifold $Y$ should be the mirror corresponding to the large complex structure limit. Otherwise what we could expect is that there exists a~{\it homological mirror}~$W$ of~$X$
equipped with a {\it non-commutative Calabi--Yau f\/ibration} $\mathrm{D^bCoh}(\PP^1)\rightarrow \mathrm{D^b}(W)$. This can be thought of as homological mirror to the Tyurin degeneration (see \cite[Section~6]{DHT} for more detials).

\section{SYZ mirror symmetry} \label{SYZ}
The Strominger--Yau--Zaslow (SYZ) mirror symmetry conjecture \cite{SYZ} provides a foundational geometric understanding of mirror symmetry for Calabi--Yau manifolds. It claims that a mirror pair of Calabi--Yau manifolds should admit dual special Lagrangian torus f\/ibrations. It is Hitchin~\cite{Hit} who f\/irst observed that the base of the f\/ibration, which is locally the moduli space of the special Lagrangian f\/ibers~\cite{McL}, carries two natural integral af\/f\/ine structures. These integral af\/f\/ine structures are essential in SYZ mirror symmetry and appear to be more fundamental than symplectic and complex geometry~\cite{GS}. One of the integral af\/f\/ine structures will play a vital role is our proof of the DHT conjecture.

\subsection{SYZ mirror symmetry for Calabi--Yau manifolds}
Let $X$ be a Calabi--Yau $n$-fold equipped a K\"ahler form $\omega$ and holomorphic volume form $\Omega$. An $n$-dimensional real submanifold $L\subset X$ is called special Lagrangian if $\omega|_L=0$ and $\operatorname{Im}(\Omega)|_L=0$ (after suitable change of the phase of $\Omega$). The celebrated SYZ mirror symmetry conjecture \cite{SYZ} asserts that, for a mirror pair of Calabi--Yau $n$-folds $X$ and $Y$,
there exist special Lagrangian $T^n$-f\/ibrations $\phi$ and $\phi^\vee$
\begin{gather*}
\xymatrix{
X \ar[rd]_\phi & & Y \ar[ld]^{\phi^\vee} \\
 & B &
}
\end{gather*}
over the common base $B$, which are f\/iberwisely dual to each other away from singular f\/ibers. The treatment of singular f\/ibers constitutes the essential part where the quantum corrections come into the play. The SYZ conjecture not only provides a powerful way to construct a mirror manifold $Y$ out of $X$ as a f\/iberwise dual, but also explains why mirror symmetry should hold via Fourier--Mukai type transformations~\cite{Leu}. It is worth noting that a mirror manifold $Y$ depends on a choice of a special Lagrangian f\/ibration on $X$, and conjecturally this is equivalent to a~choice of a large complex structure limit, where the Gromov--Hausdorf\/f limit of the Calabi--Yau mani\-fold~$X$ with its Ricci-f\/lat metric normalized to have a f\/ixed diameter is identif\/ied with the base~$B$.

\subsection[Integral af\/f\/ine structures on $B$]{Integral af\/f\/ine structures on $\boldsymbol{B}$}
Let $\phi \colon X \rightarrow B$ be a special Lagrangian $T^n$-f\/ibration of a Calabi--Yau $n$-fold $X$. We denote by $L_b$ the f\/iber of $\phi$ at $b \in B$. The complement $B^o \subset B$ of the discriminant locus carries two natural integral af\/f\/ine structures\footnote{An integral af\/f\/ine manifold is a manifold whose local coordinate changes are given by elements in the integral af\/f\/ine transformation group $\operatorname{GL}_n(\mathbb{Z})\ltimes \mathbb{R}^n$.}, which we call symplectic and complex. They are def\/ined by $\omega$ and $\operatorname{Im}(\Omega)$ respectively as follows.

{\bf Symplectic integral af\/f\/ine structure.} Let $U \subset B^o$ be a small open neighborhood of $b \in B^o$. For a basis $\{\alpha_i\}$ of $H_1(L_b,\Z)$, the symplectic integral af\/f\/ine coordinates at $b' \in U$ are def\/ined by $x_i:=\int_{A_i}\omega$, where $A_i \in H_2(X,L_b\cup L_{b'})$ is the 2-dimensional cylinder traced out by~$\alpha_i$. This is classically known in the theory of action-angle variables.

{\bf Complex integral af\/f\/ine structure.} The construction is parallel to the above. Let $U \subset B^o$ be a small open neighborhood of $b \in B^o$.
For a basis $\{\beta_i\}$ of $H_{n-1}(L_b,\Z)$, the complex integral af\/f\/ine coordinates at $b' \in U$ are def\/ined by $x^\vee_i:=\int_{B_i}\operatorname{Im}(\Omega)$, where $B_i \in H_{n}(X,L_b\cup L_{b'})$ is the $n$-dimensional cylinder traced out by $\beta_i$.

We can easily check that the above coordinates are well-def\/ined because a generic f\/iber $L$ is a special Lagrangian and both $\omega$ and $\Omega$ are closed. These coordinates depend on a choice of a~base point $b \in B^o$ and a choice of a basis of $H_1(L_b,\Z)$ or $H_{n-1}(L_b,\Z)$. Dif\/ferent base points and dif\/ferent bases of the homology group give rise to changes of coordinates by elements in~$\operatorname{GL}_n(\mathbb{Z})\ltimes \mathbb{R}^n$.

\begin{Example}
For $a,b \in \R_+$, let $X=\C/(\Z+ \sqrt{-1} a \Z)$ be the elliptic curve equipped with a holomorphic volume form $\Omega=dz$ and a K\"ahler form $\omega=\frac{b}{a}dx\wedge dy$. Then $X\rightarrow B\simeq T^1$, $z \mapsto \operatorname{Im}(z)$ def\/ines a smooth special Lagrangian $T^1$-f\/ibration.
The symplectic and complex af\/f\/ine lengths of the base $B$ are given by $b=\int_E \omega$ and $a=\int_{0}^{a}\operatorname{Im}(\Omega)$ respectively.
We observe that the mirror elliptic curve $Y=\C/(\Z+ \sqrt{-1} b \Z)$ is obtained by switching the two af\/f\/ine structures on the base~$B$. We can complexify this picture by introducing $B$-f\/ields.
\end{Example}

On the other hand, given an integral af\/f\/ine manifold $B$ of dimension $n$, we have smooth dual $T^n$-f\/ibrations:
\begin{gather*}
\xymatrix{
TB/\Lambda \ar[rd]_\phi & & T^*B/\Lambda^* \ar[ld]^{\phi^\vee}\\
 & B &
}
\end{gather*}
Here $\Lambda=\mathbb{Z}\langle \frac{\partial}{\partial x_1},\dots,\frac{\partial}{\partial x_n}\rangle$ is a f\/iberwise lattice in the tangent bundle $TB$ generated by integral af\/f\/ine coordinates $\{x_i\}$ of $B$, and $\Lambda^*$ is the dual lattice in the cotangent bundle $T^*B$. In a~natural way, $TB/\Lambda$ and $T^*B/\Lambda^*$ are complex and symplectic manifolds respectively. In order to make them (possibly non-compact) Calabi--Yau manifolds, we need a dual integral af\/f\/ine structure on~$B$. More precisely, we need a potential $f\colon B\rightarrow \mathbb{R}$ satisfying the real Monge--Amp\'ere equation $\det(\frac{\partial^2f}{\partial x_i\partial x_j} )=C$ for a constant $C \in \R$. Then the dual integral af\/f\/ine structure is given by the Legendre transformation of the original one. For example, the symplectic and complex af\/f\/ine structures discussed above are Legendre dual to each other.
This is called semi-f\/lat mirror symmetry and serves as a local model for SYZ mirror symmetry without quantum correction~\cite{Leu}. In general it is a very hard problem to extend this picture when singular f\/ibers are present.

In this article we will begin with a symplectic manifold $X$ and construct a complex mani\-fold~$Y$. So let us take a close look at this case. Given a~special Lagrangian $T^n$-f\/ibration $\phi\colon X\rightarrow B$ of a Calabi--Yau manifold $X$, we endow $B^o$ with the symplectic integral af\/f\/ine structure. We may think of the semi-f\/lat mirror $Y^o$ of $X^o:=\phi^{-1}(B^o)$ as the space of pairs $(b,\nabla)$ where~$b \in B^o$ and~$\nabla$ is a~f\/lat $\U(1)$-connection on the trivial complex line bundle over $L_b$ up to gauge. There is a~natural map $\phi^\vee\colon Y^o\rightarrow B^o$ given by forgetting the second coordinate. With the same notation as before, the complex structure of $Y^o$ is given by the following semi-f\/lat complex coordinates
\begin{gather*}
z_i(b,\nabla):=\exp\left(-2\pi\int_{A_i}\omega\right)\operatorname{Hol}_{\nabla}(\gamma_i),
\end{gather*}
where $\operatorname{Hol}_{\nabla}(\gamma_i)$ denotes the holonomy of $\nabla$ along the path $\gamma_i$. Then we observe that the dual f\/ibration $\phi^\vee$ is locally given by the tropicalization map $(z_i) \mapsto \big({-}\frac{\log|z_i|}{2\pi}\big)_i$.

\subsection{SYZ mirror symmetry for quasi-Fano manifolds} \label{SYZqFano}
Let us now consider a quasi-Fano $n$-fold $X$ with an anti-canonical divisor $Z \in |{-}K_X|$. Observing that the complement $X \setminus Z$ carries a holomorphic $n$-form with poles along~$Z$, we think of $X \setminus Z$ as a log Calabi--Yau manifold, to which the SYZ construction can be applied. More generally, in the framework of SYZ mirror symmetry for a manifold with an ef\/fective anti-canonical divisor~\cite{Aur1, CO}, the superpotential~$W$ of a mirror Landau--Ginzburg model is obtained as the weighted count of holomorphic discs of Maslov index $\mu=2$
with boundary in a smooth f\/iber $L$ of a given special Lagrangian torus f\/ibration $\phi\colon X \rightarrow B$.
To be more explicit, the superpotential $W$ is a function on the semi-f\/lat mirror $Y^o$ given by
\begin{gather*}
W(b,\nabla):=\sum_{\substack{\beta \in \pi_2(X,L_b) \\ \mu(\beta)=2}}n_\beta z_\beta(b,\nabla),
\end{gather*}
where $z_\beta$ is def\/ined to be
\begin{gather*}
z_\beta(b,\nabla):=\exp\left(-2\pi\int_{\beta}\omega\right)\operatorname{Hol}_{\nabla}(\partial \beta)
\end{gather*}
and $n_\beta$ denotes the one-point open Gromov--Witten invariant of class $\beta \in \pi_2(X,L)$ def\/ined by the machinery of Fukaya--Oh--Ohta--Ono~\cite{FOOO}. It is not dif\/f\/icult to check that $W$ is locally a~holomorphic function on~$Y^o$.

In the following, we shall focus on the toric Fano case. Namely, we consider a toric Fano $n$-fold~$X$ equipped with a toric K\"ahler form $\omega$ and a meromorphic volume form $\Omega\!=\!\wedge_{i=1}^{n}\!\sqrt{{-}1}d\log z_i$, where $(z_i)_i$ are the standard coordinates of the open dense torus $(\C^\times)^n\subset X$. Let $Z\subset X$ be the toric boundary (Example~\ref{toric}). Then the toric moment map $\phi\colon X\rightarrow \R^n$ gives a smooth special Lagrangian $T^n$-f\/ibration $\phi\colon X \setminus Z \rightarrow B^o$, where $B:=\phi(X) \subset \R^n$ is the moment polytope and~$B^o$ is its interior\footnote{By abuse of notation, we use the same $\phi$ for the restriction of $\phi$ to $X\setminus Z$.}. By construction of the semi-f\/lat mirror, it is straightforward to check the following assertion.

\begin{Proposition}\label{semi-flat}
Let $\operatorname{Trop}\colon (\C^\times)^n\rightarrow \R^n$, $(z_i)_i\mapsto \big({-}\frac{\log|z_i|}{2\pi}\big)_i$ be the tropicalization map. Then the semi-flat mirror $Y^o$ of the complement $X\setminus Z \cong (\C^\times)^n$ is given by the polyannulus $\operatorname{Trop}^{-1}(B^o)$. Moreover the dual fibration $\phi^\vee$ is identified with the restriction $\phi^\vee=\operatorname{Trop}|_{Y^o}\colon Y^o\rightarrow B^o$.
\end{Proposition}

In the toric Fano case, we do not modify $Y^o$ further, so henceforth we simply write $Y=Y^o$. In general, there is a discriminant locus in the interior of $B$ and then the semi-f\/lat mirror $Y^o$ needs quantum corrections by the wall-crossing formulae of the superpotential~$W$.

Let us take a close look at the projective line $\PP^1$. We have a special Lagrangian $T^1$-f\/ibration $\phi\colon \PP^1 \rightarrow B=[0,\operatorname{Im}(\tau)]$ given by the moment map, where $\tau:=\sqrt{-1}\int_{\PP^1}\omega$. By Proposition \ref{semi-flat}, the mirror $Y$ of $\PP^1 \setminus \{0, \infty\} \cong \C^\times$ is given by the annulus
\begin{gather*}
Y=A_{(q,1)}:=\{ q< |z| < 1\} \subset \C,
\end{gather*}
where $q:=e^{2\pi \sqrt{-1}\tau}$. Each special Lagrangian f\/iber separates $\PP^1$ into two discs, one containing~$0$ and the other containing~$\infty$. The classes $\beta_1$ and $\beta_2$ representing these disc classes satisfy $\beta_1+\beta_2=[\PP^1]$, and hence the coordinates on $Y$ should satisfy $z_{\beta_1}z_{\beta_2}=q$. Moreover we can easily check that these are the only holomorphic discs of Maslov index $2$ and $n_{\beta_1}=n_{\beta_2}=1$.
Using $z=z_{\beta_1}$ as a new coordinate on the mirror $Y=A_{(q,1)}$, we obtain the Landau--Ginzburg superpotential
\begin{gather*}
W=z_{\beta_1}+z_{\beta_2}=z+\frac{q}{z}.
\end{gather*}
So far, we discuss only the real K\"ahler structure for simplicity, but we can easily complexify it in the above discussion. We will do this in the next section.

\begin{Remark} \label{HoriVafa}
The moment map for the toric $(S^1)^n$-action is def\/ined only up to addition of a~constant in the range $\R^n$. In other words, the only intrinsic property of the base space $B$ is its af\/f\/ine structure and an af\/f\/ine embedding $B \subset \R^n$ is a choice. For example, we may take another moment map $\phi'\colon \PP^1 \rightarrow B'=\big[\frac{-\operatorname{Im}(\tau)}{2},\frac{\operatorname{Im}(\tau)}{2}\big]$, and then the mirror Landau--Ginzburg model becomes
\begin{gather*}
W'\colon \ Y'=A_{\big(q^{\frac{1}{2}},q^{-\frac{1}{2}}\big)} \rightarrow \C, \qquad z \mapsto z+\frac{q}{z},
\end{gather*}
where
\begin{gather*}
A_{(a,b)}:=\{z \in \C \, | \, a < |z| < b\}
\end{gather*}
for positive real numbers $a<b$. Note that we have a biholomorphism $A_{(q,1)} \cong A_{\big(q^{\frac{1}{2}},q^{-\frac{1}{2}}\big)}$, which is induced by the translation of the underlying af\/f\/ine manifolds $B\cong B'$ in $\R^n$. Moreover, near the large volume limit, where $\int_{\beta}\omega \rightarrow \infty$, we may identify $Y'$ with $\C^\times=A_{(0,\infty)}$. In this way, we understand that the SYZ mirror coincides with the Hori--Vafa mirror $(\C^\times,z+\frac{q}{z})$, after suitably renormalizing the superpotential (Example~\ref{HVMirror}). We refer the reader to \cite{HV} for a physical derivation of the mirror Landau--Ginzburg model.
\end{Remark}

It is also worth noting that the superpotential $W$ is in general not known to converge if $X$ is not a~toric Fano manifold. Moreover, if a~f\/iber Lagrangian $L$ bounds a holomorphic disc of class $\beta \in \pi_2(X,L)$ of Maslov index~$0$, then the one-point open Gromov--Witten invariant $n_\beta$ depends on the f\/iber $L$ as well as the point $p \in L$ which the holomorphic discs are required to pass through. On the other hand, we often want to take a smooth anti-canonical divisor $Z$ instead of the toric boundary so that the mirror Landau--Ginzburg superpotential $W$ is proper. In that case there appears to be a discriminant locus in the interior of the base $B$ and we need quantum corrections in the above toric SYZ construction~\cite{Aur1}.

\section{Gluing Landau--Ginzburg models: elliptic curves} \label{Gluing}
Finally we are in the position to conf\/irm the DHT conjecture in the case of elliptic curves. In the language of SYZ mirror symmetry, we will {\it glue the Landau--Ginzburg models} by essentially gluing the af\/f\/ine base manifolds of special Lagrangian f\/ibrations. The inspiration comes from the fact that a Tyurin degeneration is a complex analogue of a Heegaard splitting. Another key ingredient of the proof is to construct theta functions out of the Landau--Ginzburg superpotentials based on the geometry of the conjecture.

Let us consider a Tyurin degeneration of an elliptic curve $\mathcal{X}\rightarrow \Delta$, where a generic f\/iber is a~smooth elliptic curve and the central f\/iber $\mathcal{X}_0=X_1\cup_Z X_2$ is the union of two rational curves~$X_1$ and~$X_2$ glued at two points~$Z$. We complexify the K\"ahler structure $B+\sqrt{-1} \omega$ of $\mathcal{X}$ by introducing the $B$-f\/ield $B \in H^2(\mathcal{X},\R)/2 \pi H^2(\mathcal{X},\Z)$, and def\/ine
\begin{gather*}
\tau_i:=\int_{X_i}\big(B+\sqrt{-1}\omega\big),\qquad q_i:=e^{2 \pi \sqrt{-1} \tau_i}, \qquad i=1,2.
\end{gather*}
Then the complexif\/ied K\"ahler structure of a generic elliptic f\/iber $X$ of the family $\mathcal{X}\rightarrow \Delta$ is given by
\begin{gather*}
\tau :=\tau_1+\tau_2= \int_{\mathcal{X}_0} \big(B+\sqrt{-1} \omega\big)
\end{gather*}
so that $q:=e^{2 \pi \sqrt{-1} \tau}=q_1q_2$.

Let us consider the moment maps $\phi_i\colon X_i \rightarrow B_i$ for $i=1,2$,
where the base af\/f\/ine manifolds are $B_1=[0,\operatorname{Im}(\tau_1)]$ and $B_2=[-\operatorname{Im}(\tau_2),0]$.
Then the mirror Landau--Ginzburgs of $(X_1,Z)$ and $(X_2,Z)$ are respectively given by
\begin{alignat*}{4}
& W_1\colon \quad && Y_1=A_{(|q_1|,1)} \rightarrow \C, \qquad && z_1 \mapsto {z_1}+\frac{q_1}{z_1}, & \\
& W_2\colon \quad && Y_2=A_{\big(1,|q_2|^{-1}\big)} \rightarrow \C, \qquad && z_2 \mapsto z_2+\frac{q_2}{z_2}, &
\end{alignat*}
where $q_1z_2=z_1$ (Remark~\ref{HoriVafa}). We observe that the boundary of the closure $\overline{Y_1 \cup Y_2} \subset \C^\times$ can be glued by the multiplication map $z \mapsto qz$ to form an elliptic curve $Y$, which is identif\/ied with the mirror elliptic curve $\C^\times/q^{\Z}:=\C^\times/(z \sim qz)$ of~$X$. This construction corresponds to the gluing of the boundary of the af\/f\/ine manifold $B_1 \cup B_2= [-\operatorname{Im}(\tau_2),\operatorname{Im}(\tau_1)]$ by the shift $\operatorname{Im}(\tau)$ (twisted by the $B$-f\/ield upstairs).

In order to conf\/irm the DHT conjecture, we moreover want a~map $\C^\times \rightarrow (\C \times \C)\setminus (0,0)$ which descends to a double covering $W\colon \C^\times/q^{\Z} \rightarrow \PP^1$ that locally looks like the superpotential~$W_1$ (resp.~$W_2$) over the upper (resp. lower) semisphere of the base~$\PP^1$. Unfortunately, the naive analytic continuation of $(W_1,W_2)$ over $\C^\times$ does not work. The correct answer is given by considering all the Landau--Ginzburg models of the above sort, namely for $i \in \Z$ the Landau--Ginzburg models
\begin{alignat*}{4}
& W_{2i+1}\colon \quad && Y_{2i+1}=A_{\big(|q^{-i}q_1|,|q^{-i}|\big)} \rightarrow \C, \qquad && z_{2i+1} \mapsto {z_{2i+1}}+\frac{q_{1}}{z_{2i+1}}, & \\
& W_{2i}\colon \quad && Y_{2i}=A_{\big(|q^{1-i}|,|q^{1-i}q_2^{-1}|\big)} \rightarrow \C, \qquad && z_{2i} \mapsto z_{2i}+\frac{q_{2}}{z_{2i}}, &
\end{alignat*}
where the variable $z_i$ is def\/ined inductively by $qz_{i+2}=z_i$. This means to eliminate the above arbitrary choice of a fundamental domain of the $\Z$-action on $\C^\times$ to construct the mirror elliptic curve $\C^\times/q^{\Z}$. A crucial observation is that if we consider all the even or odd superpotential~$W_i$'s at once (in the sense of Remark~\ref{product formula} below), they descend to the elliptic curve $\C^\times/q^{\Z}$ as sections of an ample line bundle as follows. Let us f\/irst consider the inf\/inite product
\begin{gather*}
W'_1(z) :=\prod_{i=1}^\infty \left(1+\frac{q_1}{z_{2i-1}^2}\right)\left(1+\frac{z_{-2i+1}^2}{q_1}\right)\\
\hphantom{W'_1(z)}{}=\prod_{i=1}^\infty \big(1+q^{2i-1}\big(q_2z^2\big)^{-1}\big)\big(1+q^{2i-1}q_2z^2\big) \\
\hphantom{W'_1(z)}{} = \prod_{k=1}^\infty \frac{1}{1-q^{2k}} \sum_{l \in \Z} q^{l^2}\big(q_2z^2\big)^l \\
\hphantom{W'_1(z)}{}= \frac{e^{\frac{\pi \sqrt{-1} \tau}{6}}}{\eta(2\tau)}\vartheta_{\frac{1}{2},0}(2\zeta -\tau_1,2\tau),
\end{gather*}
where we set $z=z_1=e^{2\pi \sqrt{-1} \zeta}$ and
\begin{gather*}
\eta(\tau):= e^{\frac{\pi \sqrt{-1} \tau}{12}}\prod_{m=1}^\infty\big(1-e^{2 \pi \sqrt{-1}\tau m}\big), \qquad
\vartheta_{a,b}(\zeta,\tau):= \sum_{n \in \Z} e^{\pi \sqrt{-1} (n+a)^2\tau} e^{2\pi \sqrt{-1} (n+a)(\zeta+b)}
\end{gather*}
are the Dedekind eta function and theta function with characteristic $(a,b) \in \R^2$ respectively.

\begin{Remark} \label{product formula}
We can think of $W'_1$ as the product of all the Landau--Ginzburg superpoten\-tial~$W_i$'s for odd $i \in \Z$ because of the formula
\begin{gather*}
\left(z_j+\frac{q_k}{z_j}\right)\left(z_{-j}+\frac{q_k}{z_{-j}}\right)=q_kq^j \left(1+\frac{q_k}{z_j^2}\right)\left(1+\frac{z_{-j}^2}{q_k}\right),
\end{gather*}
for all $j \in \Z$ and $k=1,2$.
\end{Remark}

In a similar manner, we next consider
\begin{gather*}
W'_2(z) :=\prod_{i=1}^\infty \left(1+\frac{q_2}{z_{2i}^2}\right)\left(1+\frac{z_{-2i+2}^2}{q_2}\right) \\
\hphantom{W'_2(z)}{} =\prod_{i=1}^\infty \left(1+q^{2i-1}\frac{q_1}{z^2}\right)\left(1+q^{2i-1}\frac{z^2}{q_1}\right) \\
\hphantom{W'_2(z)}{} = \prod_{k=1}^\infty \frac{1}{1-q^{2k}} \sum_{l \in \Z} q^{l^2}\left(\frac{z^2}{q_1}\right)^l\\
\hphantom{W'_2(z)}{} = \frac{e^{\frac{\pi \sqrt{-1} \tau}{6}}}{\eta(2\tau)}\vartheta_{0,0}(2\zeta -\tau_1,2\tau),
\end{gather*}
which can be thought of as the product of all the Landau--Ginzburg superpotential $W_i$'s for even $i \in \Z$. It is a classical fact that the theta functions with characteristics
\begin{gather*}
\vartheta_{\frac{1}{2},0}(2\zeta -\tau_1,2\tau), \qquad \vartheta_{0,0}(2\zeta -\tau_1,2\tau)
\end{gather*}
form a basis of the $(2)$-polarization of the elliptic curve $Y=\C^\times/q^{\Z}$. Therefore we obtain a~double covering
\begin{gather*}
W\colon \ Y=\C^\times/q^{\Z} \rightarrow \PP^1, \qquad z \mapsto [W'_1(z) : W'_2(z)].
\end{gather*}
Observing that $W$ locally looks like the superpotential $W_i$ on each piece $Y_i$, we conf\/irm that it is precisely the gluing of the two Landau--Ginzburg models argued in the DHT conjecture. This completes a proof of the conjecture in the case of elliptic curves.

\begin{Remark}
It is crucial in the above proof not to take the large volume limit but to keep track of the complex structures on the mirror annuli. In this way, we are able to naturally glue the Landau--Ginzburg models (without the heuristic cutting process discussed in the conjecture). It is also interesting to observe that the product expressions of the theta functions are the manifestation of quantum corrections, which are encoded in the Landau--Ginzburg superpotentials, in SYZ mirror symmetry.
\end{Remark}

The elliptic curves are somewhat special and our construction readily generalizes to a dege\-neration of an elliptic curve to a nodal union of $n$ rational curves forming a cycle. The superpotential of each rational curve corresponds to a theta function with an appropriate characteristic, and they span a basis of the $(n)$-polarization of the mirror elliptic curve. Note that the same result is obtained in \cite[Section~8.4]{Asp} and \cite[Section~4]{KL} from dif\/ferent perspectives. This is due to the accidental fact that a Tyurin degeneration of an elliptic curve can be thought of as a~maxi\-mal degeneration at a large complex structure limit. A main dif\/ference shows up, for example, when we consider a type II degeneration of an abelian surface which is neither a~maximal nor a toric degeneration~\cite{Kan}. However the essential mechanism of the DHT conjecture is already apparent in the case of elliptic curves: gluing the base af\/f\/ine manifolds and constructing theta functions from the Landau--Ginzburg superpotentials.

\section{Further research direction} \label{Direction}
A key idea of our proof is to glue the two dif\/ferent af\/f\/ine manifolds $B_1$ and $B_2$ along the boundaries to obtain a compact af\/f\/ine manifold $\R/\operatorname{Im}(\tau)\Z$. This idea is not new and a similar construction (doubling) was already suggested by Auroux~\cite{Aur2}. However, it is in general a very hard problem to glue together higher dimensional af\/f\/ine manifolds along the boundaries. The dif\/f\/iculty is closely related to a choice of an anti-canonical divisor $Z \in |{-}K_X|$ of a quasi-Fano manifold $X$. More precisely, we need a special Lagrangian f\/ibration $\phi\colon X \rightarrow B$ which is compatible with~$Z$, and then smoothness of $Z$ is likely to be proportional to that of the boundary~$\partial B$. On the other hand, mirror symmetry for $(X,Z)$ tends to be harder as $Z$ gets less singular because we need to trade singularities of $\partial B$ with discriminant loci of the interior of~$B$. This seems the main obstruction to generalizing our discussion in higher dimensions.

A notable feature of the DHT conjecture is that it bridges a gap between mirror symmetry for Calabi--Yau manifolds and that for Fano manifolds. The conjecture naturally suggests a~construction of a mirror Calabi--Yau manifold from a more general degeneration of a Calabi--Yau manifold. For instance, the recent Gross--Siebert program~\cite{GS} is a powerful algebro-geometric program to construct a mirror Calabi--Yau manifold from a given toric degeneration of a Calabi--Yau manifold, which serves as a maximal degeneration of a Calabi--Yau manifold. However, it is not yet clear how we should make use of the mirror of a more general degeneration limit of a Calabi--Yau manifold. An expectation is that a mirror Calabi--Yau manifold should come equipped with a~map~$W$, which could be obtained by {\it gluing} (or taking the {\it linear system} of) the mirror Landau--Ginzburg superpotentials of the irreducible components of the degeneration limit. Then the glueability condition (the {\it weight} of the superpotentials being the same) should be mirror to a Kawamata--Namikawa type log smoothablity condition. In general $W$ will not provide a polarization with the mirror manifold (see also \cite[Section~5.1]{DHT}). For example it is the case if the number of the irreducible components is less than the dimension of the Calabi--Yau manifold plus~1.

Lastly, the DHT conjecture and the above speculation can be investigated from many dif\/ferent perspectives and it would be of interest to ask, for example, what they mean in the Lagrangian Floer theory and theoretical physics.

\subsection*{Acknowledgements}
The author would like to thank Yu-Wei Fan, Andrew Harder, Hansol Hong and Siu-Cheong Lau for useful conversations on related topics. Special thanks go to the referees for their valuable comments and improvements to this article. This research was supported by the Kyoto University Hakubi Project. Part of this work was carried out during the author's stay at BIRS in the fall of~2016.

\pdfbookmark[1]{References}{ref}
\LastPageEnding

\end{document}